\documentclass{amsart}
\usepackage{amsmath}
\usepackage{amsfonts}
\usepackage{amssymb}
\newtheorem{theorem}{Theorem}         
\newtheorem{corollary}{Corollary}
\newtheorem{lemma}{Lemma}

\def\star{\raise .5ex \hbox{*}}
\def\sumstar_#1{\setbox0=\hbox{$\scriptstyle{#1}$}
\setbox2=\hbox{$\displaystyle{\sum}$}
\setbox4=\hbox{${}\star\mathsurround=0pt$}
\dimen0=.5\wd0 \advance\dimen0 by-.5\wd2
\ifdim\dimen0>0pt
\ifdim\dimen0>\wd4 \kern\wd4 \else\kern\dimen0\fi\fi
\mathop{{\sum}\star}_{\kern-\wd4 #1}}
\def\sumprime_#1{\setbox0=\hbox{$\scriptstyle{#1}$}
\setbox2=\hbox{$\displaystyle{\sum}$}
\setbox4=\hbox{${}'\mathsurround=0pt$}
\dimen0=.5\wd0 \advance\dimen0 by-.5\wd2
\ifdim\dimen0>0pt
\ifdim\dimen0>\wd4 \kern\wd4 \else\kern\dimen0\fi\fi
\mathop{{\sum}'}_{\kern-\wd4 #1}}

\def\sumprime_#1{\setbox0=\hbox{$\scriptstyle{#1}$}
\setbox2=\hbox{$\displaystyle{\sum}$}
\setbox4=\hbox{${}'\mathsurround=0pt$}
\dimen0=.5\wd0 \advance\dimen0 by-.5\wd2
\ifdim\dimen0>0pt
\ifdim\dimen0>\wd4 \kern\wd4 \else\kern\dimen0\fi\fi
\mathop{{\sum}'}_{\kern-\wd4 #1}}

\def\capS{$S(t)$\ }

\title[A note on \capS  and the zeros of the Riemann zeta-function]{A note on $S(t)$  and the zeros of the Riemann zeta-function}
\author{D. A. Goldston}
\address{Department of Mathematics, San Jose
State University, San Jose, CA 95192, USA}
\email{goldston@math.sjsu.edu}
\author{S. M. Gonek}
\address{Department of Mathematics,
University of Rochester,
Rochester, NY 14627, USA}
\email{ gonek@math.rochester.edu}

\thanks{The research of both authors was supported in part by a
National Science Foundation FRG grant (DMS 0244660).  
The first author was  also partially supported by NSF grant 
DMS 0300563 and the second author was partially supported by 
NSF grant DMS 0201457. 
The authors wish to thank the Isaac Newton Institute for its hospitality
during their work on this article, and also the American Institute of Mathematics.}
\subjclass{11M26 }
\begin{document}
\maketitle

\begin{abstract}
Let $\pi S(t)$ denote the argument of the Riemann zeta-function at the 
point $\frac12+it$.  Assuming 
the Riemann Hypothesis, we sharpen the constant in the  
best currently known bounds for $S(t)$ and for the change of $S(t)$ in intervals. We then deduce estimates for the largest  multiplicity of a zero of the zeta-function and for the largest gap between the zeros.
\end{abstract}

\section{Introduction}

We assume the Riemann Hypothesis (RH) throughout this paper.

Let $N(t)$ denote the number of zeros $\rho = \frac12 +i\gamma$
of the Riemann zeta-function with ordinates in the interval $(0, t]$. 
Then for $t \geq 2$,
\begin{equation}
N(t) = \frac{t}{2\pi}\log\frac{t}{2\pi} - \frac{t}{2\pi} + 
\frac{7}{8} + S(t) + O(\frac{1}{t})\,, \label{1.1}
\end{equation}
where, if $t$ is not the ordinate of a zero,  
$S(t)$ denotes the value of 
$\frac{1}{\pi} \arg \zeta(\frac{1}{2} + i t)$
obtained by continuous variation along the straight line segments
joining $2$, $2+it$, and $\frac12+it$, starting with the value $0$.
If $t$ is the ordinate of a zero, we set $S(t)=\frac12 \lim_{\epsilon \to 0^+} 
\{S(t+\epsilon) + S(t-\epsilon)\}$. 
It follows from \eqref{1.1} that
\begin{equation}
N(t+h) - N(t) = \frac{h}{2\pi}\log\frac{t}{2\pi} + 
S(t+h) - S(t) +  O(\frac{(1+ h^{2})}{t})  \label{1.2}
\end{equation}
for $  0 < h \leq t$.
It was proved by Littlewood~\cite{L} that on RH, 
\begin{equation} 
S(t) \ll \frac{\log t}{\log\log t}. \label{1.3}
\end{equation}
Hence,  the
number of zeros with ordinates in an interval 
$(t, t+h]$ satisfies
\begin{equation}
N(t+h) - N(t) - \frac{h}{2\pi}\log\frac{t}{2\pi}   
  \ll \frac{\log t}{\log\log t}  \,, \label{1.4}
\end{equation} 
provided  that $0 < h \leq \sqrt{t}$, say. 

The bounds in \eqref{1.3} and \eqref{1.4} have not been improved over the last eighty years.
Our goal in this note is to sharpen them slightly.
\newpage
\begin{theorem} Assume the Riemann Hypothesis.  
Let $t$ be large and $0 < h \leq \sqrt{t}$. Then we have  
\begin{equation} 
\Big| N(t + h) - N( t) 
- \frac{h}{2\pi} \log \frac{t}{2\pi}\Big|  
\le \big(\frac{1}{2}+o(1)\big)\frac{\log t}{\log\log t}. \label{1.5}
\end{equation}
\end{theorem}

In light of \eqref{1.2}, this is equivalent to  
\begin{equation}
\Big| S(t + h) - S( t) \Big|  
\le \big(\frac{1}{2}+o(1)\big)\frac{\log t}{\log\log t}  \label{1.6} 
\end{equation}   
for $0 < h \leq \sqrt{t}$.  Using this, we obtain 
 
\begin{theorem} Assume the Riemann Hypothesis. Then for $t$ 
sufficiently large we have
\begin{equation}
|S(t)| \leq (\frac{1}{2} + o(1)) \frac{\log t}{\log\log t} \,. \label{1.7}
\end{equation}
\end{theorem}

To deduce Theorem 2 from Theorem 1, we use the (unconditional) estimate of Littlewood~\cite{L}
\[
\int_{0}^{T} S(u) du \ll \log T,
\]
which implies that 
\[
\int_{t}^{t+\log^2 t} S(u) du \ll \log t\,.
\]
Therefore, for $t$ sufficiently large, there is an $h$ with 
$0\leq h \leq \log^2 t$ such that $S(t+h) \leq 1$. 
Rewriting \eqref{1.6} as 
\begin{equation}
-  \big(\frac{1}{2}+o(1)\big)\frac{\log t}{\log\log t}+ S(t + h)  \leq S( t) \leq  
  \big(\frac{1}{2}+o(1)\big)\frac{\log t}{\log\log t}  + S(t + h)\, , \label{1.8}
\end{equation}   
we obtain the upper bound for $S(t)$ from the right-hand inequality.
We obtain the lower bound by using an $h$ for which  $S(t+h) \geq -1$
and the left-hand inequality. 

The following is an almost immediate corollary of Theorem 1.

\begin{corollary} Assume the Riemann Hypothesis. Let $m(\gamma)$ 
denote the multiplicity of the zero $\frac12 +i\gamma$. 
Then if $\gamma$ is sufficiently large we have
\begin{equation}
m(\gamma) \le \big(\frac{1}{2}+o(1)\big)\, 
\frac{\log \gamma}{\log\log \gamma}\,. \label{1.9} 
\end{equation}
Moreover, if $\gamma$ and $\gamma^{\prime}$ are consecutive 
ordinates and $\gamma < \gamma^{\prime}$, then
\begin{equation}
\gamma^{\prime}  - \gamma  \le 
\frac{\pi}{\log\log \gamma}(1+o(1)). \label{1.10}
\end{equation} 
\end{corollary}

To deduce \eqref{1.9} take $t = \gamma - h/2$ in
\eqref{1.5} with
$h = o(1/\log\log \gamma)$. To deduce \eqref{1.10}, assume  
$N(t+h) - N(t) =0$  in \eqref{1.5} and solve for $h$.

\section{Proof of Theorem 1}

We begin by stating two lemmas. The first lemma is a form of the Guinand-Weil explicit formula.

\begin{lemma} 
    Let $h(s)$ be analytic in the strip
    $|Im\,s| \leq \frac{1}{2}+\epsilon$ for some $\epsilon > 0$,
    and assume $|h(s)| \ll (1+|s|)^{-(1+\delta)}$  for some $\delta>0$ 
    when  $|Re\,s| \to \infty$. Let $h( w)$ be real--valued for real $w$ 
    and set $\hat{h}(x) = \int_{-\infty}^\infty h(w)
    e^{- 2\pi i x
    w } \;dw$. Then
      \begin{multline}
     \sum_{\rho} h(\frac{\rho-1/2}{i}) = h(\frac{1}{2i}) 
     + h(-\frac{1}{2i})
     - \frac{1}{2\pi}\hat{h}(0) \log \pi\\ +  
     \frac{1}{2\pi}  \int_{-\infty}^{\infty} h(u)\, 
    Re \frac{\Gamma'}{\Gamma}(\frac14 + \frac{iu}{2}) 
        \,  d u 
     - \frac{1}{2\pi} \sum_{n=2}^{\infty} 
     \frac{\Lambda(n)}{\sqrt{n} }\bigg(
     \hat{h}(\frac{\log n}{2\pi}) +
     \hat{h}(\frac{-\log n}{2\pi})  \bigg) \; .\label{2.1}
     \end{multline}
     \end{lemma}
This is a specialization of Theorem 5.12 and in particular  equation (25.10) of ~\cite{IK}.  The conditions in \cite{IK} are that $\hat{h}$ is an infinitely differentiable function with compact support, which will be satisfied in our application below; however it is not hard to also prove the lemma with the conditions we have stated.

\begin{lemma}  Let $L$  and $\delta$ be
positive real numbers and let $w=u+iv$. There exist even 
entire functions $F_{+}(w)$ and 
$F_{-}(w)$ with the following properties.
\[ \leftline{\indent ${ i) }\quad \displaystyle    F_{-}(u) \le
 \chi_{[-L,L]} (u) 
 \le F_{+}(u) \quad \text{for all real }\ u,$}\]
 \[ \leftline{\indent ${ ii) }\quad \displaystyle    
\int_{-\infty}^\infty F_+(u)\, du \le 2L + \frac{1}{\delta}, \quad 
  \int_{-\infty}^\infty F_-(u)\, du \ge 2L - \frac{1}{\delta},$}\] 
 \[ \leftline{\indent ${ iii) }\quad 
   \displaystyle  F_\pm(w) \ll e^{2\pi \delta |Im \, w|},$}\] 
\[ \leftline{\indent ${ iv) }\quad 
\displaystyle    F_{\pm}(u)\ll \min(1, \,\delta^{-2}(|u|-L)^{-2}) 
\quad \text{for}\ |u| > L,$}\] 
\[ \leftline{\indent ${ v) }\quad \displaystyle    
 \hat{F}_\pm(x)=0 \ \text{for}\ |x|\ge \delta ,$}\]
 \[ \leftline{\indent ${ vi) }\quad 
\displaystyle  \hat{F}_{\pm}(x) = \frac{\sin 2 \pi L x}{\pi x} 
+ O(\frac{1}{\delta}).$}\] 
\end{lemma}

This is essentially Lemma 2 from \cite{MonOd}. 
Functions of this type were constructed by A. Selberg, who gives a nice discussion of them in \cite{S}. For a proof of this lemma, see  H. L. Montgomery~\cite{M} and J. D. Vaaler~\cite{V}. The slightly less familiar property $iv)$  is obtained from Lemma 5 of \cite{V}.

To prove Theorem 1 we use Lemma 1 with 
$h(w) = F(w - t)$, where $t$ is large and positive
and $F$ denotes either the function  $F_{+}$ or $F_{-}$ from Lemma 2.
We assume that the parameters $\delta$ and $L$ implicit in the
definition of $F_{\pm}$ satisfy the 
conditions  
\begin{equation}\delta \geq 1, \quad \text{and} \quad 0 < L \leq 2\sqrt{t}. \label{2.2}\end{equation}
Clearly 
$ \hat{h}(x) =  e^{- 2\pi ixt} \,\hat{F}(x)$. 
Therefore, by Lemma 2 $i$) and $ii$), or $vi$),
\[
\hat{h}(0) = \hat{F}(0) = 2L + O(\frac{1}{\delta})\,,
\]
and
\[
\hat{h}( \frac{\log n}{2\pi}) 
=  n^{-it} \,\hat{F}(\frac{\log n}{2\pi}),
 \qquad
\hat{h}(-\frac{\log n}{2\pi}) 
=  n^{it} \,\hat{F}(-\frac{\log n}{2\pi}) \,.
\]
We also see that
\[
h( \frac{1}{2i})  +  h( -\frac{1}{2i}) = F(\frac{1}{2i} -t) 
+ F(\frac{1}{2i} -t)   \ll  e^{\pi \delta}.  
\]
by Lemma 2 $iii$).

We will now show that
\begin{equation}  \frac{1}{2\pi}  \int_{-\infty}^{\infty} F(u-t)\, 
    Re \frac{\Gamma'}{\Gamma}(\frac14 + \frac{iu}{2}) 
        \,  d u  = \frac{1}{2\pi} \Big(\log\frac{t}{2}\Big) \hat{F}(0) + O(1).
    \label{2.3}
\end{equation}
First, since $ Re \frac{\Gamma'}{\Gamma}(\frac14 + \frac{iu}{2})\ll \log(|u|+2)$, 
we see by \textit{iv}) of Lemma 2 and \eqref{2.2} that
\[ \begin{split} \int_{t+4\sqrt{t}}^\infty F(u-t)\, 
    Re \frac{\Gamma'}{\Gamma}(\frac14 + \frac{iu}{2}) 
        \,  d u &\ll \int_{t+4\sqrt{t}}^\infty \frac{ \log(u+2)}{\delta^2(u-t- 2\sqrt{t})^2} 
        \,  d u \\& \ll \frac{\log t}{\sqrt{t}},\end{split} \]
and similarly for the integral over $(-\infty, t-4\sqrt{t}]$. 
Next, by Stirling's formula for large $t$ together with Lemma 2 \textit{ii}) and the previous argument using \textit{iv})
\[ \begin{split} \int_{t-4\sqrt{t}}^{t+4\sqrt{t}}F(u-t)\, 
    Re \frac{\Gamma'}{\Gamma}(\frac14 + \frac{iu}{2}) 
        \,  d u&=\int_{t-4\sqrt{t}}^{t+4\sqrt{t}}F(u-t) \left( \log \frac{u}{2} + O\left(\frac{1}{1+u^2}\right) \right)\, du \\& = \int_{t-4\sqrt{t}}^{t+4\sqrt{t}}F(u-t) \left( \log \frac{t}{2} + O\left(\frac{1}{\sqrt{t}}\right) \right)\, du \\& = \Big(\log\frac{t}{2}\Big) \int_{-\infty}^{\infty}F(u-t)\, du + O\left(\frac{\log\frac{t}{2}+L}{\sqrt{t}}\right) \\ & = \Big(\log\frac{t}{2}\Big)\hat{F}(0) + O\left(1\right) .\end{split}\] 
On combining these estimates \eqref{2.3} follows.

Inserting these results into \eqref{2.1}, we obtain
\begin{equation}
\begin{aligned} 
 \sum_{\gamma} F(\gamma - t)   
=    \frac{\hat{F}(0)}{2\pi}  &\log(\frac{t}{2\pi}) 
 + O(e^{\pi\delta})  \\
&- \frac{1}{2\pi} \sum_{n=2}^\infty 
\frac{\Lambda(n)}{\sqrt{n}}  
\bigg( n^{-it} \hat{F}(\frac{\log n}{2\pi}) +  n^{it}
\hat{F}(\frac{-\log n}{2\pi})  \bigg)
 \;.
\end{aligned}\label{2.4}
\end{equation}
By $v$) and $vi$) of Lemma 2, the sum on the right is
\begin{equation*}
\begin{split} 
\ll \sum_{n < e^{2\pi\delta}} \frac{\Lambda(n)}{\sqrt{n}}  
\cos (t \log n) 
\bigg( \frac{\sin(L\log n)}{\log n} + O(\frac{1}{\delta}) \bigg) \ll \sum_{n \leq e^{2\pi\delta}} \frac{1}{\sqrt{n}}
 \ll   e^{\pi\delta} 
 \; ,
\end{split}
\end{equation*}
where the last sum was estimated trivially.
Hence,
\begin{equation}
\sum_{\gamma} F(\gamma - t)
= \frac{\hat{F}(0)}{2\pi} \log(\frac{t}{2\pi}) 
+ O(e^{\pi\delta}) \;. \label{2.5}
\end{equation}

Taking $F$ to be $F_{+}$ and using \textit{i}) and \textit{ii}) of Lemma 2,
we find that
\begin{equation*} 
\begin{split} 
N(t + L) - N(t - L)  \leq \frac{1}{2\pi}
 \log (\frac{t}{2\pi}) \bigg(2 L + \frac{1}{\delta}\bigg)
+ O(e^{\pi\delta})  \;.
\end{split}
\end{equation*}
We now take
$\pi\delta = \log \log t -2\log \log \log t$
and obtain  
\begin{equation*} 
\begin{aligned} 
N(t+ L) - 
N(t - L) \,- \,
\frac{L}{\pi}
 \log (\frac{t}{2\pi}) 
\; \leq  \; \bigg(\frac12 + o(1)\bigg) 
 \frac{\log t}{\log \log t} \,.
\end{aligned}
\end{equation*}
Had we used $F_{-}$ in \eqref{2.5} instead of $F_{+}$,
we would have found that
\begin{equation*} 
\begin{aligned} 
N(t+ L) - 
N(t - L) \,- \,
\frac{L}{\pi}
 \log (\frac{t}{2\pi}) 
\; \geq  \; \bigg(-\frac12 + o(1)\bigg) 
 \frac{\log t}{\log \log t} \,.
\end{aligned}
\end{equation*}
Combining these two inequalities, we conclude that
\begin{equation*}\label{genl upper bound}
\bigg| N(t + L) - N(t - L) \,- \,
\frac{L}{\pi} \log (\frac{t}{2\pi})  \bigg|
\leq  \; ( \frac12 + o(1) ) \,\frac{\log t}{\log \log t} \,.  
\end{equation*}
Finally, replacing $t$ by $t+ h/2$ and  taking $L = h/2$,
we obtain Theorem 1.


%

\end{document}